\documentclass[11pt,amsfonts]{article}
\usepackage{graphicx}
\usepackage{latexsym}
\usepackage{amssymb}
\usepackage{amsmath}
\usepackage{enumerate}

\usepackage{amsmath}
\usepackage{mathtools}
\usepackage{stmaryrd}
\usepackage{hyperref}
\usepackage{amssymb}
\usepackage{CJK}
\usepackage{color}

\usepackage{amsmath,latexsym,amssymb,amsthm,array,amsfonts}

\usepackage{mathrsfs,dsfont,bbm}

\usepackage{layout}

\usepackage[utf8]{inputenc}
\usepackage[LGR, T1]{fontenc}
\usepackage[english]{babel}
\usepackage{lmodern}
\usepackage{amsmath, esint}
\usepackage{amssymb}
\usepackage{color}
\usepackage{soul}
\usepackage[normalem]{ulem}
\usepackage{enumitem}
\usepackage{graphicx}
\usepackage[usenames, dvipsnames]{xcolor}
\usepackage{stmaryrd}
\usepackage{subfig}
\usepackage{multicol}
\usepackage{multirow}
\usepackage{hyperref}
\usepackage{ifthen}
\usepackage{listings}
\usepackage{diagbox}

\usepackage{tikz}

\newtheorem{prop}{Proposition}
\newtheorem{lemma}{Lemma}

\newtheorem{theorem}{Theorem}

\newtheorem{remark}{Remark}

\def\real{{\mathord{{\rm I\kern-2.8pt R}}}}        
\def\inte{{\mathord{{\rm I\kern-2.8pt N}}}}

\def\sZZ{{\rm Z\kern-2.8ptem{}Z}}

\def\z{{\mathchoice
		{\sZZ}
		{\sZZ}
		{\rm Z\kern-0.30em{}Z}
		{\rm Z\kern-0.25em{}Z} }}
\def\sQQ{{\kern 0.27em \vrule height1.45ex width0.03em depth0em
		\kern-0.30em \rm Q}}
\def\qu{{\mathchoice
		{\sQQ}
		{\sQQ}
		{\kern 0.225em \vrule height1.05ex width0.025em depth0em \kern-0.25em \rm Q}
		{\kern 0.180em \vrule height0.78ex width0.020em depth0em \kern-0.20em \rm Q}
}}
\def\sCC{{\kern 0.27em \vrule height1.45ex width0.03em depth0em
		\kern-0.30em \rm C}}
\def\complex{{\mathchoice
		{\sCC}
		{\sCC}
		{\kern 0.225em \vrule height1.05ex width0.025em depth0em \kern-0.25em \rm C}
		{\kern 0.180em \vrule height0.78ex width0.020em depth0em \kern-0.20em \rm C}
}}

\newcommand{\R}{\mathbb{R}}

\newcommand{\E}{\mathbf{E}}
\renewcommand{\P}{\mathbf{P}}

\newcommand{\norm}[2]{{\left\Vert #1 \right\Vert}_{#2}}

\newcommand{\ignore}[1]{}
\newcommand{\CQFD}{\hfill $\Box$}

\usepackage{geometry}
\begin{document}
	
	\renewcommand{\thefootnote}{\fnsymbol{footnote}}
	
	\renewcommand{\thefootnote}{\fnsymbol{footnote}}

	\title{The spatial average of solutions to  SPDEs  is  asymptotically independent of the solution}

	\author{Ciprian A. Tudor and  J\'er\'emy Zurcher\vspace*{0.2in} \\
	CNRS, Universit\'e de Lille \\
	Laboratoire Paul Painlev\'e UMR 8524\\
	F-59655 Villeneuve d'Ascq, France.\\
	\quad  ciprian.tudor@univ-lille.fr\\
	\quad jeremy.zurcher@univ-lille.fr\\
	\vspace*{0.1in} }
	
	\maketitle
	
\begin{abstract}
Let $ \left( u(t,x), t\geq 0, x\in \mathbb{R} ^{d}\right)$ be the solution to the stochastic heat or wave equation driven by a Gaussian noise which is white in time and white or correlated with respect to the spatial variable. We consider the spatial average of the solution $F_{R}(t)= \frac{1}{ \sigma _{R}}\int_{ \vert x\vert \leq R} \left( u(t,x)-1\right) dx, $
where $\sigma ^{2}_{R}= \E \left( \int_{ \vert x\vert \leq R} \left( u(t,x)-1\right) dx\right) ^{2}.$ It is known that, when $R$ goes to infinity, $F_{R}(t)$ converges in law to a standard Gaussian random variable $Z$. We show that the spatial average $F_{R}(t)$ is actually asymptotic independent by the solution itself, at any time and at any point in space, meaning that the random vector $(F_{R}(t), u(t, x_{0}))$ converges in distribution, as $R\to \infty$, to $(Z,  u(t, x_{0}))$, where $Z$ is a standard normal random variable independent of $u(t, x_{0})$.   By using the Stein-Malliavin calculus, we also obtain the rate of convergence, under the Wasserstein distance, for this limit theorem. 

\end{abstract}

	\vskip0.3cm

{\bf 2010 AMS Classification Numbers:}   60H15, 60H07, 60G15, 60F05.

\vskip0.3cm

{\bf Key Words and Phrases}: stochastic heat equation, stochastic wave equation; Malliavin derivative, Stein's method; Wasserstein distance. 
	
\section{Introduction} In the last decade, a large number of research papers analyzed the asymptotic behavior of the solution to various stochastic partial differential equations (SPDEs in the sequel).  The context is as follows. Consider a general SPDE of the form 
\begin{equation}
	\label{intro-1}
	\mathcal{L}u(t,x) =\sigma (u(t,x))\dot{ W}(t,x), \hskip0.5cm t\geq 0, x\in \mathbb{R} ^{d},
\end{equation}
with   initial condition $u(0, x)=1$ for every $x \in \mathbb{R} ^{d}$. In (\ref{intro-1}), $\mathcal{L}$ is first or second order differential operator with constant coefficients, $\sigma: \mathbb{R} \to \mathbb{R} $ is a globally Lipschitz function and $W$ is a space-time white noise (presented in the next section). The mild solution to (\ref{intro-1}) is expressed via the integral equation 
\begin{equation}\label{intro-2}
	u(t,x)= 1+ \int_{0}^{t}  \int_{\mathbb{R} ^{d}} G(t-s, x-y) \sigma (u(s,y)) W(\mathrm{d}s, \mathrm{d}y),
\end{equation}
where $G$ is the Green kernel associated with the operator $\mathcal{L}$ (i.e. the  solution of $\mathcal{L}G = \delta_0$ in the sense of distributions) and the stochastic integral with respect to $W$ is the Dalang-Walsh stochastic integral  (see Section \ref{sec21}).

Assume that the exists a unique solution to (\ref{intro-2}) (which happens for the situations considered in our work). The purpose is to study the asymptotic behavior, as $ R\to \infty$, of the quantity 
\begin{equation}
	\label{intro-3}
	F_{R}(t)= \frac{1}{ \sigma _{R}}\int_{ \vert x\vert \leq R} \left( u(t,x)-1\right) dx,
\end{equation}
where $\sigma ^{2}_{R}= \E \left( \int_{ \vert x\vert \leq R} \left( u(t,x)-1\right) dx\right) ^{2} $ and $ \vert \cdot \vert $ denotes the Euclidean norm in $ \mathbb{R} ^{d}$.  The functional $ F_{R}(t)$ is called {\it the spatial average } of the solution to the SPDE (\ref{intro-1}) at time $t$.  In  many situations, it has been proven that, for fixed $t> 0$ (in the sequel, we denote by $ \to ^{(d)}$ the convergence in distribution),
\begin{equation*}
	F_{R}(t) \to ^{(d)} Z \sim \mathcal{N}(0,1),
\end{equation*}
and  estimates for the total variation distance between the laws of $ F_{R}(t)$ and $Z\sim \mathcal{N}(0,1)$ has been obtained.  For the case of the stochastic heat equation (i.e. $\mathcal{L}=\frac{ \partial}{\partial t} -\frac{1}{2} \Delta $ in (\ref{intro-1}), with $\Delta$ the Laplacian on $\mathbb{R} ^{d}$), we refer to the papers \cite{NuHeat18}, \cite{NuHeat19} or \cite{NXZ}, which treated the situations when the noise is white or colored in space. The case of the stochastic wave equation (i.e. one takes $\mathcal{L}= \frac{ \partial ^{2}}{\partial t ^{2}}-\Delta$ in (\ref{intro-1}) and we assume in addition $ \frac{ \partial u} {\partial t} (0,x)=0$ for every $x\in \mathbb{R} ^{d}$) with white noise in time, and white or correlated spatial covariance, has been studied in  \cite{BNZ}, \cite{Delgado}, \cite{NZ}.  More precise statements of the results obtained for the heat and wave equations with white noise in time are presented in Theorem \ref{thmnuheat} and Theorem \ref{thmnuwave}.     A large variety of other situations has been considered in the literature: the case of the fractional noise in time in \cite{B1}, \cite{NSZ}, \cite{NXZ} or \cite{NZ2}, the case of the fractional heat equation (i.e. the standard Laplacian is replaced by the fractional Laplacian) in \cite{ANTV}, the case of time-independent noise in \cite{B2}, \cite{B3} and \cite{B4} or the case of delta-initial condition in \cite{CKNP}, \cite{KNP}.

In this work, we aim at studying the asymptotic independence between the spatial average $ F_{R}(t)$ given by (\ref{intro-3}) and the solution (\ref{intro-2}) itself at time $t$ and at a fixed spatial point $ x_{0} \in \mathbb{R} ^{d}$. In other words, we show that, for every $t> 0$ and $x\in \mathbb{R} ^{d}$, the family of random vectors 
\begin{equation*}
	(F_{R}(t), u(t, x_{0}), R>0)
\end{equation*}
converges in distribution, as $ R\to \infty$, to the random vector  
\begin{equation*}
	(Z, u(t, x_{0})),
\end{equation*}
where $Z\sim \mathcal{N}(0,1)$ is independent of $ u(t, x_{0})$. This implies  that the spatial average $F_{R}(t)$ is \\ asymptotically independent of the solution at any time and at any point in space. Intuitively, that means that large values of the solution are unlikely to be accompaned by large values of the spatial average, the exteme values of $u$ have no  effect on the  behavior of 
the spatial average when the integration domain in (\ref{intro-3}) becomes larger and larger. For a recent work on the concept of asymptotic independence, see \cite{DN}.

We treat several situations: the case of the (nonlinear) stochastic heat and wave equations with space-time white noise and with white noise in time and correlated spatial covariance given by the Riesz kernel.  In all these situations, we prove that the functional $ F_{R}(t)$ given by (\ref{intro-3}) is asymptotically independent of $u(t, x_{0})$ and in addition, we give estimates for the Wasserstein distance between the probability distributions of $ (F_{R}(t), u(t, x_{0}))$ and $(Z, u(t, x_{0}))$ (where $ Z \sim \mathcal{N}(0,1)$ is independent of $u(t, x_{0})$), when $R$ is large enough. 

To prove these results, we use the variant of the Stein's method recently developed in \cite{Pi} and \cite{Tud23}.  This method allows to give bounds, in terms of the Malliavin operators, for the Wasserstein distance between the law of a random vector $(X, Y)$ (with components differentiable in the Malliavin sense) and the law of the vector $(Z, Y)$, with independent components and with $Z$ a standard normal random variable. We actually slightly adapt the method from \cite{Pi} and \cite{Tud23}, by considering the situation when $X$ is a Skorohod integral. 

We organized our paper as follows. Section 2 constitutes a preliminary section where we describe  the main objects that appear in our work and we present our  new results. In Section 3 we develop our mathematical tools: we prove a new Stein-Malliavin bound to evaluate the Wasserstein distance between the law of a random vector $(X, \mathbb{Y})$, where $X$ is a Skorohod integral and $Y$ is an arbitrary $d$-dimensional vector with components differentiable in the Malliavin sense, and the vector $(Z, Y)$, where $ Z\sim \mathcal{N}(0, \sigma ^{2})$ and $Z, \mathbb{Y}$ are independent. Section 4 contains the proofs of our main results while Section 5 is the Appendix where we included the basic elements of Malliavin calculus.

\section{Preliminaries and goals of the paper}\label{sec2}

In this preliminary part, we present the stochastic heat and wave equations and some properties of their solutions. We then  recall  some results obtained in the literature regarding the asymptotic behavior of the spatial average for these equations and we  also describe our new findings.

\subsection{The Gaussian noise and the Dalang-Walsh stochastic integral}\label{sec21}

Let $(\Omega, \mathcal{F}, \P)$ a probability space. We consider an integer $d \geqslant 1$ and let $\mathcal{B}_{\text{b}} (\R^d)$ be the set of bounded Borel subsets of $ \mathbb{R} ^{d}$. 
Consider    the centered Gaussian field $W = (W(t, A), t \geqslant 0, A \in \mathcal{B}_{\text{b}} (\R^d))$ which is white in time and colored in space by the Riesz kernel. That is,  for every $s,t \geqslant 0$ and $A,B \in \mathcal{B}_{\text{b}} (\R^d)$, we have the following expression for the covariance of $W$ 
\begin{equation}\label{covbruit}
    \E [W(t, A) W(s, B)] = (s \wedge t) \int_A \int_B f_{\alpha} (x-y) \ \mathrm{d}x \ \mathrm{d}y,
\end{equation}
where the function  $f_{\alpha} : \R^d \longrightarrow \R$ is given by 
\begin{equation*}
    f_{\alpha} (x) = \left\{ \begin{array}{rl}
        \frac{1}{|x|^{d - \alpha}} & \text{ if } 0 < \alpha < d  \\
        \delta_0 (x) & \text{ if } \alpha = 0. 
    \end{array}  \right.
\end{equation*}
When $\alpha = 0$, we say that the noise is white in space, and then the covariance (\ref{covbruit}) reduces to 
\begin{equation}
    \E [W(t, A) W(s, A)] = (s \wedge t) \ \text{Leb}(A \cap B),
\end{equation}
where $\text{Leb}$ stands for the Lebesgue measure on $\R^d$. 

We define the Hilbert space $H$ as the completion of the set $\left\{ \mathbf{1}_{[0,t]} \otimes \mathbf{1}_A, t \geqslant 0, A \in \mathcal{B}_{\mathrm{b}} (\R^d) \right\}$ under the scalar product, for all $t,s \geqslant 0$ and $A,B \in \mathcal{B}_{\text{b}} (\R^d)$:
\begin{equation*}
    \left\langle \mathbf{1}_{[0,t]} \otimes \mathbf{1}_A, \mathbf{1}_{[0,s]} \otimes \mathbf{1}_B \right\rangle_H := (t \wedge s) \ \int_A \int_B f_{\alpha} (x-y) \ \mathrm{d}x \ \mathrm{d}y.
\end{equation*}
Note that if $\alpha = 0$ (the white-noise case), then $H$ is simply $L^2(\R_+ \times \R^d)$. We have, for $f,g \in H$, 
\begin{equation}\label{prodsca}
	\langle f, g\rangle _{H} = \int_{0} ^{\infty} \mathrm{d} s \int_{ \mathbb{R}^{d}} \mathrm{d}y  \int_{ \mathbb{R}^{d}} \mathrm{d}y' f(u,y) g(v, y') f_{\alpha} (y-y'). 
\end{equation}

If $f= \mathbf{1}_{[0,t]} \otimes \mathbf{1}_A,$ we set  $W(f) := W(t, A)$ for $t \geqslant 0$ and $A \in \mathcal{B}_{\mathrm{b}} (\R^d)$. Then we extend this relation to $H$, by linearity and density. In this way,  the family $(W(g), g\in H)$ becomes an isonormal process (see also \cite{T2}). We denote  by $ \mathrm{D}$ and $ \delta$ the Malliavin derivative and the Skorohod integral with respect to $W$, see Section \ref{sec:app}.

For $t>0$, denote by $ \mathcal{F}_{t}$ the $\sigma$-algebra generated by the random variables $(W(s,A), 0\leq s\leq t, A \in  \mathcal{B}_{\text{b}} (\R^d))$.  A stochastic field $(X(t,y), t\geq 0, y\in \mathbb{R} ^{d})$ is said to be $(\mathcal{F}_{t})_{t \geqslant 0}$-adapted if $X(t,x)$ is $\mathcal{F}_{t}$-measurable for every $t\geq 0$ and for every $x\in \mathbb{R}^{d}$.   Consider a random field $(X(t,y), t\geq 0, y\in \mathbb{R} ^{d}$), jointly measurable and $(\mathcal{F}_{t})_{t \geqslant 0}$-adapted  such that 
\begin{equation}\label{inte}
	\E \Vert X \Vert _{H}^{2}= \E  \int_{0} ^{\infty} \mathrm{d} u \int_{ \mathbb{R}^{d}} \mathrm{d}y  \int_{ \mathbb{R}^{d}} \mathrm{d}y' X(s,y) X(s, y')f_{\alpha}(y-y') <\infty. 
\end{equation}
For such a stochastic field, we can define the Dalang-Walsh integral (see \cite{Da} or \cite{Walsh})
\begin{equation*}
	\int_{0} ^{\infty} \int_{\mathbb{R} ^{d}} X(s,y) W (\mathrm{d}s, \mathrm{d}y).
\end{equation*}
This stochastic integral satisfies the It\^o-type isometry
\begin{equation*}
	\E \left( 	\int_{0} ^{\infty} \int_{\mathbb{R} ^{d}} X(s,y) W (\mathrm{d}s, \mathrm{d}y)\right) ^{2}= \E \Vert X \Vert _{H} ^{2},
\end{equation*}
where $\Vert \cdot \Vert _{H}= \langle \cdot, \cdot \rangle _{H}$ is the norm in $H$.

\subsection{The stochastic heat and wave equations: Convergence of  the spatial average}

This part contains some basic facts concerning the heat and wave equations. We also recall the results for the normal convergence of the spatial average, which constitute the starting of our investigation.

\subsubsection{The heat equation}

We consider the following heat equation:  for all $t \geqslant 0$ and $x \in \R^d$ :
\begin{equation}
    \frac{\partial u}{\partial t} (t, x) = \frac{1}{2} \Delta u (t, x) + \sigma(u(t, x)) \dot{W}(t,x),
\end{equation}
with initial condition $u(0, x) = 1$ for all $x \in \R^d$, where $W$ is the noise having the correlation (\ref{covbruit}), and $\sigma : \R \longrightarrow \R$ is supposed to be globally Lipschitz. We understand the solution in the \textit{mild sense}, meaning that 
\begin{equation}\label{mild11}
    u(t, x) = 1 + \int_0^t \int_{\R^d} G(t-s, x-y) \sigma(u(s, y)) \ W(\mathrm{d}s, \mathrm{d}y),
\end{equation}
where the integral must be understood in the Dalang-Walsh sense (see Section \ref{sec21})), and $G$ is the Green kernel given by, for all $t > 0$ and $x \in \R^d$,
\begin{equation}
    G(t, x) := \frac{e^{\frac{-|x|^2}{2t}}}{(2 \pi t)^{\frac{d}{2}}}.
\end{equation}
This mild solution (\ref{mild11}) exists if and only if $d < 2 + \alpha$ if $\alpha \in (0, d)$ i.e. $W$ is a white-colored noise and if and only if $d=1$  when $W$ is a space-time white noise (see for instance \cite{Da}). Note that $G$ satisfies the  semi-group property 
\begin{equation}\label{semigroup}
    \int_{\R ^{d}} G(t, x'-y) G(s, y-x) \ \mathrm{d}y = G(t+s, x' - x),
\end{equation}
for every $s, t\geq 0$ and $x, x'\in \mathbb{R}^{d}$.  We have the following bound for all $p \geqslant 1$ and $T > 0$ (see e.g. \cite{Da} or \cite{T2})
\begin{equation}
    \sup_{0 \leqslant t \leqslant T} \sup_{x \in \R^d} \E \left[ \left| u(t, x) \right|^p \right] < \infty.
\end{equation}

The object of interest of this paper is the normalized \textit{spatial average of the solution}, defined by (\ref{intro-3}).  In the case of the stochastic heat equation, we use the following notation:  for every $R > 0$ and $t \geqslant 0$, let
\begin{equation}\label{avgheat}
    F^{\text{H}}_R (t) := \frac{1}{\sigma_R} \int_{|x| \leqslant R} (u(t, x) - 1) \ \mathrm{d}x,
\end{equation}
where $\sigma_R^2 := \text{Var}\left( \int_{|x| \leqslant R} (u(t, x) - 1) \ \mathrm{d}x \right)$.  By a Fubini argument, and since the Dalang-Walsh integral coincides with the Skorohod integral for adapted integrands (see (\ref{25f-1})), we have 
\begin{equation*}
F^{\mathrm{H}}_R (t) = \delta(v_{R, t}),
\end{equation*}
 where $v_{R, t} \in L^2 (\R_+ \times \R)$ is defined as
\begin{equation}\label{vheat}
	v_{R, t} (s, y) := \frac{1}{\sigma_R} \left( \int_{|x| \leqslant R} G(t-s, x-y) \ \mathrm{d}x \right) \sigma(u(s, y)) \mathbf{1}_{[0, t]} (s).
\end{equation}

\noindent The following theorem has been proved in \cite{NuHeat18} and \cite{NuHeat19}. For now, and in  the rest of a paper, $C$ is a strictly positive constant that could change from every line to an other.

\begin{theorem}\label{thmnuheat}
    Let $\beta := d - \alpha \in (0, d]$, $T > 0$ and let $v$ be defined on (\ref{vheat}). Then for all $t > 0$ and $R > 0$, $F_R^{\mathrm{H}} (t) \in \mathbb{D}^{1,2}$. Moreover, 
\begin{enumerate}
    \item when $\beta= d=1$ (i.e. $\alpha=0$, the white-noise case),   there exists a constant $C>0$ depending on $T$ such that for all $R \geqslant 0$ and $t \in (0, T]$ 
    \begin{equation}
        d_{\mathrm{TV}} \left( F^{\mathrm{H}}_R (t), \mathcal{N}(0, 1) \right) \leqslant \sqrt{ \mathrm{Var} \left[ \left\langle \mathrm{D} F^{\mathrm{H}}_R (t), v_{R, t} \right\rangle_{L^2(\R_+ \times \R)} \right] } \leqslant \frac{C}{\sqrt{R}},
    \end{equation}
    where $d_{\mathrm{TV}}$ is the total variation distance. Moreover, as $R$ goes to $\infty$, $\sigma_R \sim C \sqrt{R}$, where $C$ depends on $T$ and $\alpha$. 

    \item for $0 < \alpha < d$ (the white-colored noise case), we have for  some $C>0$ depending on $\alpha$ and $T$ :
    \begin{equation}
        d_{\mathrm{TV}} \left( F_R^{\mathrm{H}} (t), \mathcal{N}(0, 1) \right) \leqslant \sqrt{ \mathrm{Var} \left[ \left\langle \mathrm{D} F^{\mathrm{H}}_R (t), v_{R, t} \right\rangle_{H} \right] } \leqslant \frac{C}{R^{\frac{\beta}{2}}}.
    \end{equation}
    Moreover, as $R$ goes to $\infty$, $\sigma_R \sim C R^{1 - \frac{\beta}{2}}$.
\end{enumerate}
    
\end{theorem}

For the definition of the total variation distance $d_{\mathrm{TV}}$, we refer to e.g. \cite{NuHeat18}. However, we do not need it in  our work. One last very useful result is the following lemma allowing to estimate the Malliavin derivative of $u(t, x)$ in terms of the  corresponding Green kernel. Those results are also derived in \cite{NuHeat18} and \cite{NuHeat19}.

\begin{lemma}\label{lemDnoyauheat}
    Let $T > 0$, $t \in (0, T]$ and $x \in \R^d$. Then for all $p \geqslant 2$, $u(t, x) \in \mathbb{D}^{1,p}$ and for all $s \in (0, T]$ and $y \in \R^d$:
    \begin{equation}\label{Dnoyauheat}
        \E \left[ \left| \mathrm{D}_{s,y} u (t, x) \right|^p \right]^{\frac{1}{p}} \leqslant C_{p, T} G(t-s, x-y).
    \end{equation}
\end{lemma}   

\subsubsection{The wave equation in dimension 1}

We consider the following one-dimensional wave equation: 
\begin{equation}
    \frac{\partial^2 u}{\partial t^2} (t, x) = \frac{\partial^2 u}{\partial x^2} (t, x) + \sigma(u(t, x)) \dot{W}(t,x), \hskip0.4cm t\geq 0, x\in \mathbb{R},
\end{equation}
with initial condition $u(0, x) = 1$ and $\frac{\partial u}{\partial t}(0, x) = 0$ for all $x \in \R$, where $W$ is the noise having the covariance (\ref{covbruit}), and $\sigma : \R \longrightarrow \R$ is supposed to be globally Lipschitz. We understand again the solution in the \textit{mild sense}, meaning that 
\begin{equation}
    u(t, x) = 1 + \int_0^t \int_{\R^d} G_1 (t-s, x-y) \sigma(u(s, y)) \ W(\mathrm{d}s, \mathrm{d}y),
\end{equation}
where the integral  $W(\mathrm{d}s, \mathrm{d}y)$ is the Dalang-Walsh integral introduced in Section \ref{sec21} and $G_1$ is the Green kernel given by for all $t > 0$ and $x \in \R$
\begin{equation}\label{Gwave}
    G_1 (t, x) := \frac{1}{2} \mathbf{1}_{\{ |x| \leqslant t \}}.
\end{equation}
Then we have the following bound for all $p \geqslant 1$ and $T > 0$ (see \cite{Da} or \cite{T2})
\begin{equation}\label{uLp}
    \sup_{0 \leqslant t \leqslant T} \sup_{x \in \R} \E \Big[ \big| u(t, x) \big|^p \Big] < \infty.
\end{equation}
For every $R > 0$ and $t \geqslant 0$, we consider the  spatial average of the solution 
\begin{equation}\label{avgwave}
    F_R^{\text{W}} (t) := \frac{1}{\sigma_R} \int_{-R}^{R}  (u(t, x) - 1) \ \mathrm{d}x,
\end{equation}
where $\sigma_R^2 := \text{Var}\left( \int_{-R}^R (u(t, x) - 1) \ \mathrm{d}x \right)$. We can write  \begin{equation*}
	F_R^{\mathrm{W}} (t) = \delta(v_{R, t}),
\end{equation*}
 where 
\begin{equation}\label{vwave}
	v_{R, t} (s, y) := \frac{1}{\sigma_R} \left( \int_{-R}^R G_1 (t-s, x-y) \ \mathrm{d}x \right) \ \sigma(u(s, y)) \ \mathbf{1}_{[0, t]} (s). 
\end{equation}

The following theorem has been proved in the reference \cite{Delgado}.
\begin{theorem}\label{thmnuwave}
    Let $\beta := 1 - \alpha \in [0, 1)$, and $v$ defined on (\ref{vwave}).
\begin{enumerate}
    \item When $\beta = 1$ (the white-noise case), we have for all $t> 0$ and $R$ large enough 
    \begin{equation}
        d_{\mathrm{TV}} \left( F_R^{\mathrm{W}} (t), \mathcal{N}(0, 1) \right) \leqslant \sqrt{\mathrm{Var} \left[ \left\langle \mathrm{D} F^{\mathrm{W}}_R (t), v_{R, t} \right\rangle_{L^2(\R_+ \times \R)} \right]} \leqslant \frac{C}{\sqrt{R}},
    \end{equation}
    Moreover, as $R$ goes to $\infty$, $\sigma_R \sim C \sqrt{R}$. 

    \item For $0 < \alpha < 1$ (the white-colored noise case), we have for all $t> 0$ and $R$ large enough 
    \begin{equation}
        d_{\mathrm{TV}} \left( F^{\mathrm{W}}_R (t), \mathcal{N}(0, 1) \right) \leqslant \sqrt{\mathrm{Var} \Big[ \left\langle \mathrm{D} F_R^{\mathrm{W}} (t), v_{R, t} \right\rangle_{H} \Big]} \leqslant \frac{C}{R^{\frac{\beta}{2}}}.
    \end{equation}
    Moreover, as $R$ goes to $\infty$, $\sigma_R \sim C R^{1 - \frac{\beta}{2}}$.
\end{enumerate}
    
\end{theorem}

As for the heat equation we have a lemma giving an estimate of the  Malliavin derivative of the solution (it has been proven in \cite{Delgado}, Lemma 2.2). 

\begin{lemma}\label{lemDnoyauwave}
    Let $T > 0$, $t \in (0, T]$ and $x \in \R^d$. Then for all $p \geqslant 2$, $u(t, x) \in \mathbb{D}^{1,p}$ and for all $s \in (0, T]$ and $y \in \R^d$ :
    \begin{equation}\label{Dnoyauwave}
        \E \Big[ \big| \mathrm{D}_{s,y} u (t, x) \big|^p \Big]^{\frac{1}{p}} \leqslant C_{p, T} G_1 (t-s, x-y).
    \end{equation}
\end{lemma}   

\subsection{Results}

We define first the distance we will work with. For every random variable $X$, real or vector valued, we will denote by $\P_X$ or $\mathcal{L}(X)$ its law. We define $\text{Lip}(1)$ the set of functions $h : \R \times \R^d \longrightarrow \R$ which are Lipschitz continuous, with a Lipschitz constant less or equal to $1$. In other words,
\begin{eqnarray*}
    \text{Lip}(1) := \left\{ h : \R \times \R^d \longrightarrow \R \left| \sup_{\substack{\mathbf{x}, \mathbf{x}' \in \R^{d+1} \\ \mathbf{x} \neq \mathbf{x}'}} \frac{|h(\mathbf{x}) - h(\mathbf{x}')|}{|\mathbf{x} - \mathbf{x}'|} \leqslant 1 \right. \right\}.
\end{eqnarray*}
We define the \textit{Wasserstein distance} for every random variable $\mathbb{X},\mathbb{Y}$ valued in $\R \times \R^d$ such that $h(\mathbb{X}) \in L^1 (\P)$ and $h(\mathbb{Y}) \in L^1 (\P)$ as :
$$ d_{\mathrm{W}} \big( \mathcal{L}(\mathbb{X}), \mathcal{L}(\mathbb{Y}) \big) := \sup_{h \in \text{Lip}(1)} \Big| \E[h(\mathbb{X})] - \E[h(\mathbb{Y})] \Big|. $$

We state the first main result which shows that the spatial average of the solution of the heat equation is asymptotically independent of the solution itself. 

\begin{theorem}\label{thmheat}
    Let $T > 0$, $F^{\mathrm{H}}_R = \left( F^{\mathrm{H}}_R (t), t \in [0, T] \right)$ given by (\ref{avgheat}) and let  $x_0 \in \R$ be an arbitrary point. $\beta := d - \alpha$.
    \begin{enumerate}
        \item If $\beta = d= 1$, we have for every $t \in (0, T]$ and for $R$ large,
        \begin{equation}\label{dWheatblanc}
            d_{\mathrm{W}} \Big( \mathcal{L} \left( F^{\mathrm{H}}_R(t), u(t, x_0) \right), \mathcal{N}(0, 1) \otimes \mathcal{L}(u(t, x_0)) \Big) \leqslant \frac{C}{\sqrt{R}},
        \end{equation}
        where $C>0$ depends on $T$.
        \item If $0 < \alpha < d$, then for $R$ large enough,
        \begin{equation}\label{dWheatcolor}
            d_{\mathrm{W}} \left( \mathcal{L} \left( F^{\mathrm{H}}_R(t), u(t, x_0) \right), \mathcal{N}(0, 1) \otimes \mathcal{L}(u(t, x_0)) \right) \leqslant \frac{C}{R^{\frac{\beta}{2}}},
        \end{equation}
    where  $C>0$ depends on $T$ and $\alpha$.
    \end{enumerate}
\end{theorem}
 A similar phenomenon happens in the case of the stochastic wave equation in spatial dimension $d=1$.

\begin{theorem}\label{thmwave}
    Let $T > 0$, $F^{\mathrm{W}}_R = \left( F^{\mathrm{W}}_R (t), t \geqslant 0 \right)$ given by (\ref{avgwave}) and $x_0 \in \R$.  Let $\beta := 1 - \alpha$
    \begin{enumerate}
        \item If $\beta = 1$, \textit{i.e.} in  the  space-time white noise case, there exists a constant $C > 0$ depending on $T$ and $x_0$ such that for every $t \in (0, T]$ and for all $R > x_0 + 2T$ :
        \begin{equation}\label{dWwaveblanc}
            d_{\mathrm{W}} \left( \mathcal{L} \left( F^{\mathrm{W}}_R (t), u(t, x_0) \right), \mathcal{N}(0, 1) \otimes \mathcal{L}(u(t, x_0)) \right) \leqslant \frac{C}{\sqrt{R}};
        \end{equation}
        \item If $0 < \alpha < 1$ (the white colored noise case), then there exists a constant $C > 0$ depending on $T$, $\alpha$ and $x_0$ such that for every $t \in (0, T]$ and for all $R > x_0 + 2T$ :
        \begin{equation}\label{dWwavecolor}
            d_{\mathrm{W}} \left( \mathcal{L} \left( F^{\mathrm{W}}_R (t), u(t, x_0) \right), \mathcal{N}(0, 1) \otimes \mathcal{L}(u(t, x_0)) \right) \leqslant \frac{C}{R^{\frac{\beta}{2}}}.
        \end{equation}
   
    \end{enumerate}
\end{theorem}

\section{A fundamental Stein's bound}

We begin our proofs by introducing a new Stein's bound for a couple of random variables, by slightly adapting a result from \cite{Tud23}. The main result of this section is Proposition \ref{steinbound}, where we give a bound for the Wasserstein distance between the law of $(X, \mathbb{Y})$, where is $X$ is a Skorohod integral, and the law of $(Z, Y)$ with $ Z\sim \mathcal{N} (0, \sigma ^{2})$ independent of $\mathbb{Y}$. 
The idea to prove it is to use the following  observation  (see for instance Lemmas 1 and 2 in \cite{Tud23},  see also \cite{Pi}): if $X$ is a real random variable and $\mathbb{Y}$ a random vector, then the couple $(X, \mathbb{Y})$ is independent and $X$ follows $\mathcal{N}(0, \sigma^2)$ if and only if for every function $f : \R \times \R^d \longrightarrow \R$ belonging to a large class of function, we have 
\begin{equation*}
    \sigma^2 \E \left[ \frac{\partial f}{\partial x}(X, \mathbb{Y}) \right] = \E [X f(X, \mathbb{Y})].
\end{equation*}
To the above identity, one can associate a Stein's equation as follows: let  $h : \R \times \R^d \longrightarrow \R$ such that $\E[|h(Z, \mathbb{Y})|] < \infty$, where $Z \sim \mathcal{N}(0, \sigma^2)$ is independent of $\mathbb{Y}$. We set,  for all $x \in \R$ and $\mathbf{y} \in \R^d$:
\begin{equation}\label{eqstein}
    \sigma^2 \frac{\partial f}{\partial x} (x, \mathbf{y}) - x f(x, \mathbf{y}) = h(x, \mathbf{y}) - \E[h(Z, \mathbf{y})]. 
\end{equation}
 We will use the following lemma, proved in Proposition 1 in \cite{Tud23}.
 
\begin{lemma}\label{lembound} Let  $h: \mathbb{R} \to \mathbb{R} ^{d}$ be of $C^1$ with partial derivatives bounded by $1$, the equation (\ref{eqstein}) admits a unique bounded solution. Moreover, there exists a constant $C > 0$ which does not depend of $h$ such that
$$ \max \left\{ \norm{\frac{\partial f_h}{\partial x}}{\infty}, \max_{1 \leqslant j \leqslant d} \norm{\frac{\partial f_h}{\partial y_j}}{\infty}  \right\} \leqslant C. $$
\end{lemma}

We can state and prove the following proposition, which constitutes a generalization of Proposition 2.2 in \cite{NuHeat18}.

\begin{prop}\label{steinbound}
    Let $X = \delta(v)$, where $v \in \mathrm{Dom}(\delta) \subset L^2 (\Omega \times H)$ and let  $\mathbb{Y} = (Y_1, \cdots, Y_d)$ be a vector-valued random variable. Suppose that $\E[X^2] = \sigma^2$, $X \in \mathbb{D}^{1,2}$ and $\mathbb{Y}=(Y_{1},...,Y_{d})$ with $ Y_{j} \in \mathbb{D}^{1,2}$ for each $j=1,...,d$.  Then we have 
    \begin{equation}\label{diststein}
        d_{\mathrm{W}} \left( \mathcal{L}(X, \mathbb{Y}), \mathcal{N} \left(0, \sigma^2 \right) \otimes \mathcal{L}(\mathbb{Y}) \right) \leqslant C \left\{ \sqrt{\mathrm{Var} \left[ \langle \mathrm{D}X, v \rangle_{H} \right]} + \sum_{j=1}^d \sqrt{\E \left[ \left\langle v, \mathrm{D} Y_j \right\rangle^2_H \right]} \right\}.
    \end{equation}
\end{prop} 

\noindent \textbf{Proof.} Suppose first that $h$ is  also $C^1$ on $\mathbb{R} ^{d+1}$ and it has partial derivatives bounded by $1$. Let $Z \sim \mathcal{N} \left( 0, \sigma^2 \right)$ be a random variable independent of $\mathbb{Y}$. Then, by the Stein's equation (\ref{eqstein})
\begin{eqnarray*}
    \E[h(X, \mathbb{Y})] - \E[h(Z, \mathbb{Y})] = & & \E \Big[ h(X, \mathbb{Y}) - \E[h(Z, \mathbb{Y})] \Big] \\
    = & & \E \left[ \sigma^2 \frac{\partial f_h}{\partial x} (X, \mathbb{Y}) - X f_h (X, \mathbb{Y}) \right] \\
    = & & \E \left[ \sigma^2 \frac{\partial f_h}{\partial x} (X, \mathbb{Y}) \right] - \E \left[ \delta(v) f_h (X, \mathbb{Y}) \right]. 
\end{eqnarray*}
By duality formula (\ref{eq: duality formula}) for the divergence operator and by chain rule of Malliavin derivative 
\begin{eqnarray*}
    \E[h(X, \mathbb{Y})] - \E[h(Z, \mathbb{Y})] = & & \E \left[ \sigma^2 \frac{\partial f_h}{\partial x} (X, \mathbb{Y}) \right] - \E \left[ \left\langle v, \mathrm{D} f_h (X, \mathbb{Y}) \right\rangle_H \right] \\
    = & & \E \left[ \frac{\partial f_h}{\partial x} (X, \mathbb{Y}) \left( \sigma^2 - \langle v, \mathrm{D}X \rangle \right) \right] - \sum_{j=1}^d \E \left[ \frac{\partial f_h}{\partial y_j} (X, \mathbb{Y})  \langle v, \mathrm{D}Y_j \rangle  \right].
\end{eqnarray*}
Using again  the  duality  relationship (\ref{eq: duality formula}),  we have 
$$\E[\langle v, \mathrm{D}X \rangle] = \E[\delta(v)^2] = \E[X^2] = \sigma^2, $$
we conclude by Cauchy--Schwarz inequality and by Lemma \ref{lembound} that 
\begin{equation}\label{inghbornee}
    \Big| \E[h(X, \mathbb{Y})] - \E[h(Z, \mathbb{Y})] \Big| \leqslant C \left\{ \sqrt{\mathrm{Var} \left[ \langle \mathrm{D}X, v \rangle_{H} \right]} + \sum_{j=1}^d \sqrt{\E \left[ \left\langle v, \mathrm{D} Y_j \right\rangle^2_H \right]} \right\},
\end{equation}
where the constant $C$ does not depend on $h$. This inequality is true for every $h \in \mathrm{Lip}(1)$ being $C^1$. For a general $h \in \mathrm{Lip}(1)$, we approximate it for every $\varepsilon > 0$ by 
\begin{eqnarray*}
    h_{\varepsilon} (x, \mathrm{y}) := \E \left[ h \left( x + \sqrt{\varepsilon} N, \mathbf{y} + \sqrt{\varepsilon} \mathbf{N} \right) \right],
\end{eqnarray*}
where $N \sim \mathcal{N}(0, 1)$ and $\mathbf{N} \sim \mathcal{N}_d (0, I_d)$ are independent. Then $h_{\varepsilon} \in \mathrm{Lip}(1)$ is $C^1$ and uniformly converges to $h$ on $\R \times \R^d$ as $\varepsilon$ goes to zero. Consequently, we get 
\begin{equation*}
    \Big| \E[h(X, \mathbb{Y})] - \E[h(Z, \mathbb{Y})] \Big| \leqslant C \left\{ \sup_{\R \times \R^d} \big| h - h_{\varepsilon} \big| + \sqrt{\mathrm{Var} \left[ \langle \mathrm{D}X, v \rangle_{H} \right]} + \sum_{j=1}^d \sqrt{\E \left[ \left\langle v, \mathrm{D} Y_j \right\rangle^2_H \right]} \right\}.
\end{equation*}
Since this holds for every $\varepsilon > 0$, we conclude on (\ref{inghbornee}) for every $h \in \mathrm{Lip}(1)$ by using the uniform convergence of $h_{\varepsilon}$  to $h$, and so in (\ref{diststein}). \CQFD

\section{Proofs}
This section is consecrated to the proofs of our main results, Theorem \ref{thmheat} and Theorem \ref{thmwave}. The two proofs follow the same main line but the calculations are pretty different since they strongly depend on the expression of the Green kernel. 

\subsection{Proof of Theorem \ref{thmheat}}
We separate the proofs of the cases when the noise is white  or correlated with respect to the space variable. 

\subsubsection{White noise case}

We suppose that $\alpha = 0$ and so $d = 1$. Let $t \in [0, T]$ and $R > 0$. Recall that  $F^{\mathrm{H}}_R (t) = \delta(v_{R, t})$ where $v_{R, t} \in L^2 (\R_+ \times \R)$ is defined by (\ref{vheat}).
Then, since $F^{\mathrm{H}}_R (t) \in \mathbb{D}^{1,2}$, we apply Proposition \ref{steinbound} to get 
\begin{eqnarray}\label{dWHeat1}
    & & d_{\mathrm{W}} \left( \mathcal{L} \left( F^{\mathrm{H}}_R(t), u(x_0, t) \right), \ \mathcal{N}(0, 1) \otimes \mathcal{L} \big( u(x_0, t) \big) \right) \\ \nonumber
    \leqslant & & C \left\{ \sqrt{\mathrm{Var} \left[ \langle \mathrm{D} F^{\mathrm{H}}_R(t), v_{R, t} \rangle_{H} \right]} + \sqrt{\E \left[ \left\langle  v_{R, t}, \mathrm{D} u(x_0, t) \right\rangle_{H}^2 \right]} \right\}.
\end{eqnarray}
By Theorem \ref{thmnuheat}, the first term satisfies
\begin{equation*}
    \sqrt{\mathrm{Var} \left[ \langle \mathrm{D} F^{\mathrm{H}}_R(t), v_{R, t} \rangle_{H} \right]} \leqslant \frac{C}{\sqrt{R}}.
\end{equation*}
All we need to do is to estimate the second term. By definition, we have
\begin{equation*}
    \Big\langle v_{R, t}, \mathrm{D} u(t, x_0) \Big\rangle_{H} = \frac{1}{\sigma_R} \int_{-R}^R \left[ \int_0^t \int_{\R} \sigma(u(s, y)) G(t-s, x-y) \ \mathrm{D}_{s,y} u(t, x_0) \ \mathrm{d}y \mathrm{d}s \right] \mathrm{d}x.
\end{equation*}
Hence, the expectation of the square of this term yields to, by using Cauchy--Schwarz inequality :
\begin{eqnarray*}
    \E \left[ \Big\langle v_{R, t}, \mathrm{D} u(t, x_0) \Big\rangle_{H} ^2 \right] = & & \frac{1}{\sigma_R^2} \int_{-R}^R \int_{-R}^R \ \int_0^t \int_{\R} \ \int_0^t \int_{\R} G(t-s, x-y) G(t-s', x'-y') \\
    & \cdot & \E \Big[ \sigma(u(s, y)) \mathrm{D}_{s, y} u (t, x_0) \ \sigma(u(s', y')) \mathrm{D}_{s', y'} u (t, x_0) \Big] \ \mathrm{d}y \mathrm{d}s \ \mathrm{d}y' \mathrm{d}s' \ \mathrm{d}x \mathrm{d}x' \\
    \leqslant & &  \frac{1}{\sigma_R^2} \int_{-R}^R \int_{-R}^R \ \int_0^t \int_{\R} \ \int_0^t \int_{\R} G(t-s, x-y) G(t-s', x'-y') \\
    & \cdot & \E \left[ \sigma(u(s, y))^4 \right]^{\frac{1}{4}} \E \left[ \sigma(u(s', y'))^4 \right]^{\frac{1}{4}} \E \left[ \mathrm{D}_{s, y} u (t, x_0)^4 \right]^{\frac{1}{4}}  \E \left[ \mathrm{D}_{s', y'} u (t, x_0)^4 \right]^{\frac{1}{4}} \\
    & & \mathrm{d}y \mathrm{d}s \ \mathrm{d}y' \mathrm{d}s' \ \mathrm{d}x \mathrm{d}x'.
\end{eqnarray*}
Since $\sigma$ is globally Lipschitz (we denote by $\norm{\sigma}{\mathrm{Lip}}$ its Lipschitz constant), we have by using (\ref{uLp}) that for all $s \in [0, t]$ and $y \in \R$:
\begin{equation}\label{bornesigma}
    \sup_{0 \leqslant t \leqslant T} \sup_{y \in \R} \E \left[ \sigma(u(s, y))^4 \right] \leqslant 8 |\sigma(0)|^4 + 8 \norm{\sigma}{\text{Lip}}^4 \sup_{0 \leqslant t \leqslant T} \sup_{x \in \R} \E \left[ u(t, x)^4 \right] < \infty.
\end{equation}
Moreover, by Lemma \ref{lemDnoyauheat}, we have for all $s \in [0, t]$ and $y \in \R$:
\begin{equation*}
    \E \left[ \mathrm{D}_{s, y} u (t, x_0)^4 \right]^{\frac{1}{4}} \leqslant C G(t-s, x_0-y).
\end{equation*}
Consequently, we have 
\begin{eqnarray*}
     \E \left[ \Big\langle v_{R, t}, \mathrm{D} u(t, x_0) \Big\rangle_{H}^2 \right] \leqslant & &  \frac{C}{\sigma_R^2} \int_{-R}^R \int_{-R}^R \ \int_0^t \int_{\R} \ \int_0^t \int_{\R}  G(t-s, x-y) G(t-s', x'-y') \\
     & &  G(t-s, x_0-y) G(t-s', x_0-y') \ \mathrm{d}y \mathrm{d}s \ \mathrm{d}y' \mathrm{d}s' \ \mathrm{d}x \mathrm{d}x'.
\end{eqnarray*}
By semi-group property (\ref{semigroup}):
\begin{equation*}
     \E \left[ \Big\langle v_{R, t}, \mathrm{D} u(t, x_0) \Big\rangle_{H}^2 \right] \leqslant  \frac{C}{\sigma_R^2} \int_{-R}^R \int_{-R}^R \ \int_0^t \ \int_0^t  G(2(t-s), x-x_0) G(2(t-s'), x'-x_0) \ \mathrm{d}s \mathrm{d}s' \ \mathrm{d}x \mathrm{d}x'.
\end{equation*}
To conclude, we do the majoration $\int_{-R}^{R} \leqslant \int_{\R}$, we compute $\int_{\R} G(2(t-s), x) \ \mathrm{d}x = 2 (t-s)$ and we find
\begin{equation*}
     \E \left[ \Big\langle v_{R, t}, \mathrm{D} u(t, x_0) \Big\rangle_{H}^2 \right] \leqslant  \frac{C}{\sigma_R^2} \leqslant \frac{C}{R},
\end{equation*}
by Theorem \ref{thmnuheat}. Then the estimate  (\ref{dWheatblanc}) is obtained.

\subsubsection{White-colored  noise case}

We suppose now that $d \geqslant 1$ and $0 < \alpha < d$. We still have $F_R^{\mathrm{H}} (t) = \delta(v_{R, t})$, with $v$ defined in (\ref{vheat}), and the estimation (\ref{dWHeat1}). By Theorem \ref{thmnuheat}, we have 

\begin{equation*}
    \sqrt{\mathrm{Var} \left[ \langle \mathrm{D} F^{\mathrm{H}}_R(t), v_{R, t} \rangle_H \right]} \leqslant \frac{C}{R^{\frac{\beta}{2}}}.
\end{equation*}
So we focus on the second term. By the definition of the scalar product on $H$ (see (\ref{prodsca})),
\begin{equation*}
    \Big\langle v_{R, t}, \mathrm{D} u(t, x_0) \Big\rangle_H = \frac{1}{\sigma_R} \int_{|x| \leqslant R} \left[ \int_0^t \int_{\R^d} \int_{\R^d} \sigma(u(s, y)) G(t-s, x-y) \ \mathrm{D}_{s,z} u(t, x_0) \ \frac{\mathrm{d}y \ \mathrm{d}z}{|y-z|^{\beta}} \ \mathrm{d}s \right] \mathrm{d}x.
\end{equation*}
By using Cauchy--Schwarz inequality, (\ref{bornesigma}) and Lemma \ref{lemDnoyauheat}, we have
\begin{eqnarray*}
&&    \E \left[ \Big\langle v_{R, t}, \mathrm{D} u(t, x_0) \Big\rangle_H^2 \right]\\
 && \leqslant  \frac{C}{\sigma_R^2} \int_{|x| \leqslant R} \int_{|x'| \leqslant R} \ \int_0^t \int_{\R^d} \int_{\R^d} \  \int_0^t \int_{\R^d} \int_{\R^d} G(t-s, x-y) G(t-s', x'-y') \\ 
    & & G(t-s, x_0 - z) G(t-s', x_0 - z') \frac{\mathrm{d}y \ \mathrm{d}z}{|y-z|^{\beta}} \ \mathrm{d}s \frac{\mathrm{d}y' \ \mathrm{d}z'}{|y'-z'|^{\beta}} \ \mathrm{d}s' \ \mathrm{d}x' \ \mathrm{d}x \\
     &&= \frac{C}{\sigma_R^2} \left( \int_{|x| \leqslant R} \ \int_0^t \int_{\R^d} \int_{\R^d}  G(t-s, y)   G(t-s, z) \frac{\mathrm{d}y \ \mathrm{d}z}{|z-y + (x-x_0)|^{\beta}} \ \mathrm{d}s  \ \mathrm{d}x \right)^2 \\
     &&= \frac{C}{\sigma_R^2} \left( \int_{|x| \leqslant R} \ \int_0^t \E \left[  \frac{1}{ \left| x-x_0-\sqrt{t-s}(Z_1 - Z_2) \right|^{\beta}}  \right] \ \mathrm{d}s  \ \mathrm{d}x \right)^2,
\end{eqnarray*}
where $Z_1$ and $Z_2$ are independent and follows the standard $d$-dimensional Gaussian law $\mathcal{N}_d (0, I_d)$. By doing the change of variables $x \leftarrow \frac{x}{R}$ and using the fact that
\begin{equation*}
    \sup_{y \in \R^d} \int_{|x| \leqslant 1} \frac{\mathrm{d}x}{|x+y|^{\beta}} < \infty,
\end{equation*}
we conclude that 
\begin{eqnarray*}
    \E \left[ \Big\langle v_{R, t}, \mathrm{D} u(t, x_0) \Big\rangle_H^2 \right] \leqslant & & \frac{C R^{2d - 2 \beta}}{\sigma_R^2} \left( \int_0^t \E \left[ \int_{|x| \leqslant 1}  \frac{\mathrm{d}x}{ \left| x-\frac{x_0 + \sqrt{t-s}(Z_1 - Z_2)}{R} \right|^{\beta}} \right] \ \mathrm{d}s   \right)^2 \\
    \leqslant & &  \frac{C R^{2d - 2 \beta}}{\sigma_R^2}.
\end{eqnarray*}
By Theorem \ref{thmnuheat}, $\sigma_R^2$ is equivalent to $C R^{2d-\beta}$ for $R$ large, so we conclude that
\begin{equation*}
     \E \left[ \Big\langle v_{R, t}, \mathrm{D} u(t, x_0) \Big\rangle_H^2 \right] \leqslant \frac{C}{R^{\beta}}.
\end{equation*}
This leads to (\ref{dWheatcolor}) and concludes the proof.

\subsection{Proof of Theorem \ref{thmwave}}

The idea of the proof is the same as for  the previous theorem, except that the semi-group property does not hold for the Green kernel associated to the  wave equation. Instead, we will compute explicitly the integrals appearing in the proof. 

\subsubsection{White noise case}

Suppose that $\alpha = 0$ and $d=1$. Then we write $F_R^{\mathrm{W}} (t) = \delta(v_{R, t})$, where 
$v_{R, t} $ is given by (\ref{vwave}).  Then, by doing the same steps as in the proof of Theorem \ref{thmheat}, by using Theorem \ref{thmnuwave} and Proposition \ref{steinbound}, we only need to focus on estimating $\E \left[ \left\langle v_{R, t}, \mathrm{D} u(t, x_0) \right\rangle^2 \right]$. With the use of Cauchy--Schwarz inequality, inequality (\ref{bornesigma}) that still holds for the wave equation and Lemma \ref{lemDnoyauwave}, we conclude that
\begin{eqnarray*}
    \E \left[ \Big\langle v_{R, t}, \mathrm{D} u(t, x_0) \Big\rangle^2 \right] \leqslant & & \frac{C}{\sigma_R^2} \int_{-R}^R \int_{-R}^R \ \int_0^t \int_{\R} \ \int_0^t \int_{\R} G_1(t-s, x-y) G_1(t-s', x'-y') \\
    & & G_1(t-s, x_0-y) G_1(t-s', x_0-y') \ \mathrm{d}y \mathrm{d}s \ \mathrm{d}y' \mathrm{d}s' \ \mathrm{d}x \mathrm{d}x' \\
    = & & \frac{C}{\sigma_R^2} \left( \int_{-R}^R \ \int_0^t \int_{\R}  G_1(t-s, x-y) G_1(t-s, x_0-y) \ \mathrm{d}y \mathrm{d}s \ \mathrm{d}x \right)^2. 
\end{eqnarray*}
We go on by computing this integral, by using the explicit expression of $G_1$ given in (\ref{Gwave}). Then, we have for $s \in [0, t]$ and $x \in \R$ that 
\begin{eqnarray*}
    & & \int_{\R} G_1 (t-s, x-y) G_1 (t-s, x_0-y) \ \mathrm{d}y \\
    = & &  \Big( x_0-x-2(t-s) \Big) \ \mathbf{1}_{[x_0 - 2(t-s), x_0]} (x) + \Big( x-x_0-2(t-s) \Big) \ \mathbf{1}_{[x_0, x_0  + 2(t-s)]} (x),
\end{eqnarray*}
where, by convention, $[a,b] = \varnothing$ if $a > b$. Hence, we have for $R > x_0 + 2t$ : 
\begin{eqnarray*}
     \E \left[ \Big\langle v_{R, t}, \mathrm{D} u(t, x_0) \Big\rangle^2 \right] \leqslant & & \frac{C}{\sigma_R^2} \left[ \int_0^t \int_{\max\{ -R, x_0 - 2(t-s) \}}^{\min\{ x_0, R \}} (x_0-x-2(t-s)) \ \mathrm{d}x \mathrm{d}s \right. \\
    & + & \left. \int_0^t \int_{\max\{ -R, x_0 \}}^{\min\{ x_0 + 2(t-s), R \}} (x-x_0-2(t-s)) \ \mathrm{d}x \mathrm{d}s  \right] \\
    = & &  \frac{C}{\sigma_R^2} \left[ \int_0^t \int_{ x_0 - 2(t-s)}^{x_0} (x_0-x-2(t-s)) \ \mathrm{d}x \ \mathrm{d}s \right. \\
    & + & \left. \int_0^t \int_{x_0}^{x_0 + 2(t-s)} (x-x_0-2(t-s)) \ \mathrm{d}x \ \mathrm{d}s  \right] \\
    = & & \frac{C}{\sigma_R^2}.
\end{eqnarray*}
By Theorem \ref{thmnuwave}, $\sigma_R^2$ is equivalent to $C R$ when $R$ goes to infinity, so we conclude that 
\begin{equation*}
     \E \left[ \Big\langle v_{R, t}, \mathrm{D} u(t, x_0) \Big\rangle^2 \right] \leqslant \frac{C}{\sqrt{R}},
\end{equation*}
which proves the  inequality (\ref{dWwaveblanc}).

\subsubsection{White-colored noise case}

Suppose that $0 < \alpha < 1$. By repeating the same steps as before, we have 
\begin{equation}\label{EG1G1}
    \E \left[ \Big\langle v_{R, t}, \mathrm{D} u(t, x_0) \Big\rangle_H^2 \right] \leqslant \frac{C}{\sigma_R^2} \left( \int_{-R}^R \int_0^t \int_{\R} \int_{\R} \frac{G_1 (t-s, x-y) G_1 (t-s, x_0-z)}{|y-z|^{\beta}} \ \mathrm{d}y \mathrm{d}z \ \mathrm{d}s \ \mathrm{d}x \right)^2.
\end{equation}
We need to compute the double integral with respect to $y$ and $z$.

\begin{lemma}
    Let $x \in \R$ and $s \in [0, t]$. Then
    \begin{eqnarray*}
        & & (1- \beta)(2 - \beta) \int_{\R} \int_{\R} \frac{G_1 (t-s, x-y) G_1 (t-s, x_0-z)}{|y-z|^{\beta}} \ \mathrm{d}y \mathrm{d}z \\
        = & & \left( (x_0 - x + 2(t-s))^{2-\beta} + (x_0 - x - 2(t-s))^{2-\beta} - 2 (x_0 - x)^{2 - \beta} \right) \ \mathbf{1}_{(- \infty, x_0 - 2(t-s))} (x) \\
        & + &  \left( (x_0 - x + 2(t-s))^{2-\beta} + (x - x_0 + 2(t-s))^{2-\beta} - 2 (x_0 - x)^{2 - \beta} \right) \ \mathbf{1}_{[x_0 - 2(t-s), x_0)} (x) \\
        & + &  \left( (x - x_0 + 2(t-s))^{2-\beta} + (x_0 - x + 2(t-s))^{2-\beta} - 2 (x - x_0)^{2 - \beta} \right) \ \mathbf{1}_{(x_0, x_0 + 2(t-s)]} (x) \\
        & + & \left( (x - x_0 + 2(t-s))^{2-\beta} + (x - x_0 - 2(t-s))^{2-\beta} - 2 (x - x_0)^{2 - \beta} \right) \ \mathbf{1}_{(x_0 + 2(t-s), + \infty)} (x).
    \end{eqnarray*}
\end{lemma}
\noindent Going back on (\ref{EG1G1}), we actually have with the previous lemma, with $R > x_0 + 2t$ :
\begin{eqnarray*}
    & & \E \left[ \Big\langle v_{R, t}, \mathrm{D} u(t, x_0) \Big\rangle_H^2 \right] \\
    \leqslant & & \frac{C}{\sigma_R^2} \left( C' + \int_0^t \int_{-R}^{x_0-2(t-s)} \left( (x_0 - x + 2(t-s))^{2-\beta} + (x_0 - x - 2(t-s))^{2-\beta} \right. \right. \\
    & & \left.  - 2 (x_0 - x)^{2 - \beta} \right) \ \mathrm{d}x \mathrm{d}s  \\
    & + & \left. \int_0^t \int_{x_0+2(t-s)}^R \left( (x - x_0 + 2(t-s))^{2-\beta} + (x - x_0 - 2(t-s))^{2-\beta} - 2 (x - x_0)^{2 - \beta} \right) \ \mathrm{d}x \mathrm{d}s \right)^2,
\end{eqnarray*}
where $C'$ is another constant (with respect to $R$) that may depend on $\alpha, T, x_0$. Then
\begin{eqnarray*}
    & & \E \left[ \Big\langle v_{R, t}, \mathrm{D} u(t, x_0) \Big\rangle_H^2 \right] \\
    \leqslant & & \frac{C}{\sigma_R^2} \left( C' + \int_0^t \left( (x_0 + R + 2(t-s))^{3-\beta} + (x_0 + R - 2(t-s))^{3-\beta}   - 2 (x_0 + R)^{3 - \beta} \right) \ \mathrm{d}s \right. \\
    & + & \left. \int_0^t \left( (R - x_0 + 2(t-s))^{3-\beta} + (R - x_0 - 2(t-s))^{3-\beta} - 2 (R - x_0)^{3 - \beta} \right) \ \mathrm{d}s \right)^2 \\
    = & & \frac{C}{\sigma_R^2} \Bigg\{ C' + \frac{1}{2(4-\beta)} \left( (x_0 + R + 2t)^{4 - \beta} - (x_0 + R)^{4 - \beta} - (x_0 + R - 2t)^{4 - \beta} + (x_0 + R)^{4 - \beta} \right)  \\
    & + &   \frac{1}{2(4 - \beta)} \left( (R - x_0 + 2t)^{4-\beta} - (R - x_0)^{4-\beta}  - (R - x_0 - 2t)^{4-\beta} + (R-x_0)^{4 - \beta} \right)     \\
    & + & 2t(x_0 + R)^{3 - \beta} - 2t (R - x_0)^{3 - \beta} \Bigg \}^2 \\
    = & & \frac{C}{\sigma_R^2} \Bigg\{ C' + \frac{1}{2(4-\beta)} \left( (R + x_0 + 2t)^{4-\beta} - (R - (x_0 + 2t))^{4 - \beta}  \right. \\
    & + & \left. (R + 2t - x_0)^{4 - \beta}   - (R - (2t - x_0))^{4-\beta}  \right) - 2t(x_0 + R)^{3 - \beta} - 2t (R - x_0)^{3 - \beta} \Bigg \}^2.
\end{eqnarray*}
By expanding $(R + u)^{4 - \beta} = R^{4 - \beta} + u(4-\beta) R^{3 - \beta} + u^2 \frac{(4-\beta)(3-\beta)}{2} R^{2 - \beta} + O(R^{1 - \beta})$ when $R$ grows to infinity, and expanding $(u + R)^{3 - \beta}$ too, we conclude that 
\begin{equation*}
    \E \left[ \Big\langle v_{R, t}, \mathrm{D} u(t, x_0) \Big\rangle_H^2 \right] \leqslant \frac{C}{\sigma_R^2} \left( C' + C'' R^{1 - \beta} \right) \leqslant \frac{C R^{2 - 2 \beta}}{\sigma_R^2}.
\end{equation*}
By Theorem \ref{thmnuwave}, $\sigma_R^2$ is equivalent to $C R^{2 - \beta}$ for $R$ large. Consequently, we have
\begin{equation*}
    \E \left[ \Big\langle v_{R, t}, \mathrm{D} u(t, x_0) \Big\rangle_H^2 \right] \leqslant \frac{C}{\sigma_R^2} \left( C' + C'' R^{1 - \beta} \right) \leqslant \frac{C}{R^{\beta}}.
\end{equation*}
This concludes in (\ref{dWwavecolor}), and the proof.

\vskip0.5cm 

Let us end this section with some comments:

\begin{remark}
	\begin{itemize}
		\item In Theorems \ref{thmheat} and \ref{thmwave}, the distance between the probability laws of $ (F_{R}(t), u(t, x_{0}))$ and 
	the random vector  $	(Z, u(t, x_{0})),$  where $Z\sim \mathcal{N}(0,1)$ is independent of $ u(t, x_{0})$, is the sum of two terms: the  first quantifies Wasserstein distance  between $F_{R}(t)$ and $Z$ and the second somehow measures the correlation between the spatial average $ F_{R}(t)$ and $ u(t, x_{0})$ when $R$ is larger and larger.  It can be noticed from the proofs of Theorems \ref{thmheat} and \ref{thmwave} that these two terms are of the same order when the integration domain becomes larger and larger (i.e. $R\to \infty$).
	
	\item In the case of the stochastic wave equation, we restricted to the situation when $d=1$. This is because our calculations strongly rely on the expression of the Green kernel $G_{1}$. On the other hand, it seems that the main idea used in this work could be extended to the case when the spatial dimension is bigger than one. 
	\end{itemize}
\end{remark}

\section{Appendix: Malliavin calculus}\label{sec:app}
Let us next describe the basic tools from Malliavin calculus needed in this work. We refer to \cite{N} for a detailed presentation. Let $ (W(h), h\in H)$ be an isonormal process, i.e. a centered Gaussian family such that  for every $h, g\in H$,
\begin{equation*}
	\mathbf{E} W(h) W(g)= \langle h, g\rangle_{H}, 
\end{equation*}We introduce
$C_p^{\infty}(\mathbb{R}^n)$ as the space of smooth functions with all their partial derivatives having at most polynomial growth at infinity, and $\mathcal{S}$ as the space of simple random variables of the form 
\begin{equation*}
	F = f(W(h_1), \dots, W(h_n)),
\end{equation*}
where $f\in C_p^{\infty}(\mathbb{R} ^{n})$ and $h_i \in H$, $1\leq i \leq n$. Then the Malliavin derivative $DF$ is defined as $H$-valued random variable
\begin{equation}
	DF=\sum_{i=1}^n  \frac {\partial f} {\partial x_i} (W(h_1), \dots, W(h_n)) h_i\,.
\end{equation}
For any $p\geq 1$, the operator $D$  is closable as an operator  from $L^p(\Omega)$ into $L^p(\Omega;  H)$. Then $\mathbb{D}^{1,p}$ is defined as the completion of $\mathcal{S}$ with respect to the norm
\begin{equation*}
	\|F\|_{1,p} = \left(\E [|F|^p] +   \E(  \|D F\|^{p}_{H})   \right)^{1/p}\,.
\end{equation*}
The adjoint operator $\delta$ of the derivative is defined through the duality formula
\begin{equation}\label{eq: duality formula}
	\E (\delta(u) F) = \E( \langle u, DF \rangle_{H}),
\end{equation}
valid for any $F \in \mathbb{D}^{1,2}$ and any $u\in {\rm Dom} \, \delta \subset L^2(\Omega; H) $. The operator $\delta$ is
also called the Skorokhod integral since, in the case of the standard Brownian motion,
it coincides
with an extension of the It\^o integral introduced by Skorokhod (see e,g, \cite{N}).

In our context, the Hilbert space $H$ coincides with $ L ^{2}(\mathbb{R}_{+}\times \mathbb{R})$ when $W$ is space-time white noise and it is characterized by the  scalar product (\ref{prodsca}) when $W$ is a white noise in time with spatial covariance given by the Riesz kernel.  In particular, for any   adapted random field $X$ which is jointly measurable and satisfies (\ref{inte}) belongs to the domain of $\delta$, and $\delta (X)$ coincides with the Walsh integral (introduced in Section \ref{sec2}).
\begin{equation}
\label{25f-1}
\delta (X) = 
\int_0^\infty \int_{\R} X(s,y) W(\mathrm{d} s, \mathrm{d} y).
\end{equation}

\noindent {\bf Funding: } Both authors  acknowledge partial support from  Labex CEMPI (ANR-11-LABX-007-01). C. Tudor also acknowledges support from   the  ANR project SDAIM 22-CE40-0015,  ECOS SUD (project C2107), Japan Science and Technology Agency CREST  (grant JPMJCR2115) and  by the Ministry of Research, Innovation and Digitalization (Romania), grant CF-194-PNRR-III-C9-2023.

\end{document}